\documentclass[12pt,reqno,a4paper]{amsart}

\usepackage{euscript,fullpage,amsmath,amssymb,mathabx,enumerate,xcolor,float}
\usepackage[T1]{fontenc}

\newcommand\numberthis{\addtocounter{equation}{1}\tag{\theequation}}

\newtheorem{thm}{Theorem}[section]

\newtheorem{remark}[thm]{Remark}

\newtheorem{proposition}[thm]{Proposition}

\newtheorem{definition}{Definition}
\newtheorem*{standing}{Standing Assumptions}

\newcommand{\LL}{\mathcal{L}}
\newcommand{\LLL}{\tilde{\mathcal{L}}_n}

\newcommand{\al}{\alpha}

\DeclareMathOperator{\essinf}{ess\ inf}

\title{Targets and holes}
\date{\today}
\begin{document}

\begin{abstract}
We address the extreme value problem of a one-dimensional dynamical system approaching a fixed target while constrained to avoid a fixed set---which can be thought of as a small hole. The presence of the latter influences the extremal index which will now depend explicitly on the escape rate.
\end{abstract}

\keywords{Extreme value theory, Hitting time statistics, Open systems, Escape rates}
\subjclass[2010]{Primary: 37A25, 60G70; Secondary: 37E05, 37D50}

\maketitle
\begin{center}
\authors{P. Giulietti\footnote{Scuola Normale Superiore - Centro di Ricerca Matematica Ennio De Giorgi, Piazza dei
Cavalieri 7, 56126 Pisa, Italy. E-mail: {\tt
\email{paolo.giulietti@sns.it}}.},
P. Koltai \footnote{Department of Mathematics and Computer Science, Freie
Universit\"at Berlin, Arnimallee 6, 14195 Berlin, Germany. E-mail: {\tt
\email{peter.koltai@fu-berlin.de. }}},
S.\ Vaienti\footnote{Aix Marseille Universit\'e, Universit\'e de Toulon, CNRS, CPT, 13009 Marseille, France. E-mail: {\tt \email{vaienti@cpt.univ-mrs.fr}}.}}

\end{center}

%\tableofcontents
\section{Introduction}

This work is motivated by the appearance of extreme events in specific natural contexts. We are interested in the statistical description of phenomena where a perishable dynamics (i.e., an open system) is approaching a fixed target state.
As examples one can think of the process describing a hurricane approaching a city or a pandemic outbreak (with the underlying space being the spatial distribution) approaching a critical extension, before they disappear.
Thus, the dynamical setting is novel in that it has two main features: in the phase space, on one hand there is a target point which will be approximated by small balls around it, and on the other hand there is an absorbing region which terminates the process on entering it.

A one dimensional prototype of such situation can be formulated as an \emph{extreme value problem for an open system}, thus allowing a rigorous study. Similar setups, restricted to the presence of  shrinking targets or absorbing regions, but not both, have already been studied in many situations; see~\cite{BDT,DT,FFT,LM} for a comprehensive account of the literature.

We consider a dynamical system where there is an absorbing region, a \emph{hole} $H$, such that an orbit entering terminates its evolution (i.e., it is lost forever).
By considering the orbits of the whole state space, it is possible to construct a \emph{surviving set}.  On this we fix a point and a small ball around it, the \emph{target set}~$B$. We investigate the probability of hitting $B$ for the first time after $n$ steps while avoiding~$H$, in the $n\to\infty$ limit. We will show that this question can be formulated in a precise probabilistic manner by introducing conditionally invariant probability measures for the open system.

For the purpose of our motivation, we will call the entrance of the system trajectory into the target an extreme event, and the closest approach of the trajectory to the target is measured by so-called extreme values (of a suitable function of the distance).
An extreme value distribution (EVD) will be obtained by means of a spectral approach on suitably perturbed transfer
operators (see, among others, \cite{S1, S2, S3, S4, S5, GK}). The boundary levels and the extremal index of the EVD will be expressed in terms of the Hausdorff dimension of the surviving set and of the escape rate, respectively.
The EVD will explicitly depend on whether the target point in the surviving set is periodic or not, cf.\ our main result, Proposition~\ref{prop:final}.
The theory above can also be adapted to handle a sequence of
target sets which shrink to a point outside the surviving set. In this case, it predicts correctly that the EVD is degenerate, i.e., the dynamics cannot approach the target point indefinitely. These three cases together thus define a \emph{trichotomy} of possible EVDs.

%\pk{Is this ok?:}
Our approach also links parameters of the EVD to dynamical quantities, thus it provides tools of computing dynamical indicators through approximating the limiting distribution by the so-called Generalized Extreme Value (GEV) distribution, and vice a versa; this will be the object of future investigations.

In Section~\ref{sec:framework} we will detail the systems we will consider.
Section \ref{sec:evd} will present the deduction of the extreme value distribution by using a well-established spectral approach. In Section~\ref{sec:ei} we will compute explicitly the extremal index. The full statement of the result is Proposition~\ref{prop:final} in Section~\ref{ssec:PTP}. Last, in Section~\ref{hf} we show how a degenerate EVD arises
when the target set becomes disjoint from the surviving sets.
For the sake of simplicity, we will restrict ourselves to uniform expanding maps of the intervals, although generalizations are possible following the same approach. The remarks after Proposition~\ref{prop:final} discuss possible extensions.

\section{The open system} \label{sec:framework}
To access open systems through an operator-theoretic framework, we will adapt the theory developed by C.~Liverani and V.~Maume-Deschamps~\cite{LM}. They considered Lasota--Yorke maps\footnote{I.e., uniformly expanding maps, $\inf_I |T'|=\beta>1$, such that there exists a finite partition of the interval $I$ with the property that $T$ restricted to the closure of each element is $C^1$ and monotone.} $T:I \righttoleftarrow $ on the unit interval $I$ and a transfer operator with a potential $g$ of bounded variation (BV).

We denote with $\LL_g$  the transfer (Perron--Frobenius) operator associated to
$T$ and $g$; it acts on functions $f \in BV \cap L^1(\mu_g) $
as 
\begin{equation} \label{transferop}
\LL f(x)=\sum_{Ty=x}f(y)g(y),
\end{equation}
where $\mu_g$ is the conformal measure left invariant by the  dual $\LL^*$ of the transfer operator,
\[
\LL^* \mu_g = e^{P(g)} \mu_g,
\]
where $P(g)$ is the topological pressure of the potential~$g$.

For simplicity, we will restrict ourselves to the potential $g=\frac{1}{\left|T'\right|}$, however we refer to
Remarks~\ref{rem:extension1}-\ref{rem:extension2} below for possible extensions of the result to general potentials. First of all note that, in this case, the conformal measure $\mu_{\left|T'\right|^{-1}}$ will be Lebesgue (denoted by $m$) and~$P(g)=0$. Recall that if we equip the space of $\text{BV}$ functions with the norm given by the total
variation plus the $L^1$ norm\footnote{From now on we will denote
will denote $L^{1}(m)$ and  $L^{\infty}(m)$ by  $L^1$ and $L^{\infty}$. The
$L^1$ norm will be written as $|\cdot|_1.$}, then the unit ball of such a
$\text{BV}$ space  is compact in $L^1$; this will allow us to make good use of the
spectral decomposition of transfer operators. Moreover, our
probability distributions will be explicitly written  in terms of the Lebesgue
measure and therefore they will be accessible to numerical computations.
We will use later on the quantity $\Theta(g)$ defined as $\log\Theta(g):=\lim_{n\rightarrow \infty}\frac1n \log \sup_I g_n,$ where  $g_n=g(x)\times \cdots \times  g(T^{n-1}x)$;  in our case it simply becomes $\Theta(g)=\beta.$

We then consider a proper subset $H\subset I$ of measure $0<m(H)<1$, called
the \emph{hole}, and its complementary set $X_0=I\setminus H.$ We denote by
$\smash{X_n=\bigcap_{i=0}^n T^{-i}X_0}$ the set of points that have not yet fallen into the hole at time~$n$.  The \emph{surviving set} will be denoted by~$\smash{X_{\infty} = \bigcap_{n=1}^{\infty} X_n}$. The key object in our study are conditionally invariant  probability measures.

\begin{definition} \label{def:cipm}
A probability measure $\nu$ which is absolutely continuous with respect to Lebesgue is called a conditionally invariant  probability measure if it satisfies for any Borel set $A\subset I$  and for all $n>0$ that
\begin{equation} \label{eq:cipm}
\nu(T^{-n}A\cap X_n)= \nu(A)\ \nu(X_n).
\end{equation}
We use for it the abbreviation a.c.c.i.p.m.
\end{definition}
The measure $\nu$ is supported on $X_0$, $\nu(X_0)=1$, and moreover
\[
\nu(X_n)= \alpha^n, \ \text{where}\ \nu(X_1) = \nu(T^{-1}X_0)=\alpha < 1 .
\]
Apart from being absolutely continuous with respect to Lebesgue, this measure is numerically accessible in simulations. The existence of a.c.c.i.p.m. in our setting is achieved by Theorem A in \cite{LM}.
Note that $\alpha$ contains, at the same time, the information about $m(H)$ and the expansion of the system (see equations (\ref{CA1}) and (\ref{CA2}) below).
  We now introduce our first perturbed transfer
operator defined on bounded variation function $f$ as
\begin{equation} \label{eq:elle0}
\LL_0(f)=\LL(f{\bf 1}_{X_0}).
\end{equation}

We will use the following facts which are summarized in
\cite[Lemma 1.1, Lemma 4.3]{LM}:
\begin{itemize}
\item Let $\nu={\bf 1}_{X_0} h_0m$ with $h_0\in L^1$ then $\nu$ is an a.c.c.i.p.m.\
if and only if $\LL_0 h_0= \alpha h_0,$ for some $\alpha\in (0,1]$.
\item Let $\alpha,h_0$ be as above. Moreover, let
$\mu_0$ be a probability measure on $I$ such that $\LL^{*}_0 \mu_0=\alpha \mu_0.$
Then $\mu_0$ is supported in $X_{\infty}$\footnote{This result implies that $X_{\infty}$ is not empty. Otherwise, this fact follows by compactness whenever all the $X_n$ were closed, but this is not always the case if at least one branch of $T$ is not onto. } and the
measure $\Lambda$ with
\[
\Lambda=h_0\, \mu_0 \quad \text{is } T\text{-invariant}.
\]
\item The measure $\mu_0$ satisfies the conformal relation:
\begin{equation}\label{CM}
\mu_0(TA)=\alpha\int_A  |T'|  d\mu_0,
\end{equation}
for every measurable set $A\subset I$ on which $T$ is one to one.    

\item For any $v\in L^1(\mu_0)$ and $w\in L^{\infty}(\mu_0)$ we have the duality
relationship:
       \begin{equation}\label{DU}
      \int \LL_0 v \ w \ d\mu_0=\alpha \int v \ w\circ T\ d\mu_0.
      \end{equation}
      Actually this duality formula will only be used to rewrite the integral (\ref{QQQQ}) and in that case $w$ will be the characteristic function of a measurable set and $v=h_0$ which is $\mu_0$-integrable.
\end{itemize}
We are now strengthening our assumptions by taking small holes since in this case we can use the results in \cite[Section 7]{LM} and that will allow us to  apply the spectral approach of extreme value theory.
 We first need a few
preparatory results which will be also essential for the next considerations.

\subsection{Lasota--Yorke inequalities}
\label{ssec:LSI}

Lemma 7.4 in \cite{LM}  states that for each $\chi\in (\beta=\Theta(g), 1)$
there exists $a, b>0$, independent of $H$, such that,
for each $w$ of bounded
variation:
\begin{align}\label{LY1}
\| \LL^nw \|_{BV}\le a\chi^n \|w\|_{BV}+b|w|_1\\
\| \LL_0^nw \|_{BV}\le a\chi^n \|w\|_{BV}+b|w|_1.
\end{align}
The proof of the first inequality is standard; the second one relies on the fact that the jumps in the total variation norm of the  backward images of the hole grow linearly with $n$ and they are dominated by the exponential contraction of the derivative, see also the proof of~\cite[Theorem~2.1]{BV}.

\subsection{Closeness of the transfer operators and their spectra}
\label{ssec:closeness}

We introduce a so-called triple norm, defined by
$\vvvert \mathcal{P} \vvvert_1:=\sup_{ \|w \|_{BV} \le 1}|\mathcal{P}w|_1,$ where $w\in
\text{BV}$ and the linear operator $\mathcal{P}$ maps into~$L^1.$\footnote{If we use a different measure ``$\text{meas}$'' instead of $m$ we
will write $ \vvvert\mathcal{P} \vvvert_{\text{meas}}$.} It is easily proven in~\cite[Lemma 7.2]{LM} that
\begin{equation} \label{CO}
 \vvvert \LL-\LL_0  \vvvert_1 \le e^{P(g)}m(H)=m(H).
\end{equation}
The idea is now to take a hole of small $m$-measure in such a way that even the spectra of the two operators are close. This is achieved next.

The following result is proved in~\cite[Theorem 7.3]{LM}. For each $\chi_1\in (\chi, 1)$ and $\delta\in (0,1-\chi_1),$ there exists $\epsilon_0>0$ such that if $ \vvvert\LL_0-\LL \vvvert_1\le \epsilon_0$ then the spectrum of $\LL_0$ outside the disk $\{z\in \mathbb{C}, |z|\le \chi_1\}$ is $\delta$-close, with multiplicity, to the one of $\LL.$ This result will allow us to get a very useful quasi-compactness representation for the two operators, which will be the starting point of the perturbation theory of extreme values.

\subsection{Quasi-compactness of the transfer operators}

First of all we should add a further prescription for our unperturbed system, namely we will require that $T$ has a unique invariant measure $\mu$ absolutely continuous with respect $m$ with density $h$ and moreover the system $(I, T, \mu)$ is mixing.
Therefore $\LL h= h$ and since $\LL^*m=m,$ we have that $\mu=h m.$ Moreover, for any function $v$ of bounded variation, there exists a linear operator $\mathcal{Q}$ with spectral radius $\text{sp}(\mathcal{Q})$ strictly less than $1$, such that
\begin{equation}\label{QC1}
\LL v=h\int v \, dm+\mathcal{Q}v.
\end{equation}
By the closeness of the spectra the same representation holds for $\LL_0$, namely there will be a number $\lambda_0,$ a non-negative function $h_0$ of bounded variation, a  probability measure $\mu_0$ and a linear operator $\mathcal{Q}_0$ with spectral radius strictly less than $\lambda_0$ such that for any $v\in \text{BV}$:
\begin{align}\label{QC2}
\LL_0 h_0=\lambda_0 h_0, \ \LL_0^*\mu_0=\lambda_0\mu_0\\
\lambda_0^{-1}\LL_0 v=h_0\int v \, d\mu_0+\mathcal{Q}_0v.
\end{align}
Notice that with what we discussed above, in the framework of small holes we will have $\lambda_0=\alpha$ and moreover the measure $\Lambda=h_0\mu_0$ will be $T$-invariant.

\section{Extreme Value Distribution} \label{sec:evd}
For a fixed \emph{target point} $z\in X_{\infty}$ let us consider the observable
\[
\phi(x)=-\log|x-z| \quad \text{for $x\in I$},
\]
and the function
$$
M_n (x):=\max\{\phi(x), \cdots, \phi(T^{n-1}x)\}.
$$
For $ u \in \mathbb{R}_+$, we are interested in the probabilities of $M_n\le u,$ where $M_n$ is now seen as a random variable on a suitable (yet to be chosen) probability space $(\Omega, \mathbb{P}).$ First of all we notice that the set of $x \in I$ for which it holds $\{M_n\le u \}$ is equivalent to the set  $\{ \phi\le u, \ldots, \phi\circ T^{n-1}\le u\}$. In turn this is the set $E_n:=(B^c\cap  T^{-1}B^c\cdots \cap  T^{-(n-1)} B^c)$ where, for simplicity of notation, we denote with $B^c$ the complement of the open ball $B:=B(z, e^{-u})$, which we call the \emph{target} (set). So far we are following points which will enter the ball $B$ for the first time after at least $n$ steps, but we should also guarantee that they have not fallen into the hole before entering the target. Therefore we should consider the event: $E_n\cap X_{n-1}$  conditioned on~$X_{n-1}$, i.e., conditioned on the event of not terminating at least for $n-1$ steps. To assure that, the natural sequence of probability measures is given by the following
\begin{definition}
For any Borel set $A\subset I$ and any $n\ge 1$ we introduce the sequence of  probability measures:
\[
\mathbb{P}_n(A):=\frac{\nu(A\cap X_{n-1})}{\nu(X_{n-1})}.
\]
\end{definition}

Suppose now that, rather than taking one ball $B,$ we consider a sequence of
balls $B_n:=B(z, e^{-u_n})$ centered at  the target point $z$ and
of radius $e^{-u_n}.$
 Therefore:
\begin{equation}\label{SEI}
\mathbb{P}_n(M_n\le u_n) = \frac{1}{\nu(X_{n-1})}\int_I {\bf 1}_{B_n^c\cap X_0}(x)\cdots {\bf 1}_{B_n^c\cap X_0}(T^{n-1}x)d\nu,
\end{equation}
and we will consider the limit for  $n\rightarrow \infty$, where $u_n$ is a \emph{boundary level} sequence which guarantees the existence of a non-degenerate limit. We anticipate that such a sequence will be dictated directly by the proof below and it must satisfy for a given~$\tau$
\begin{equation}\label{BBLL}
n \ \Lambda \big(B(z, e^{-u_n} ) \big)\rightarrow \tau\quad\text{as }n\to\infty.
\end{equation}
By introducing our second perturbed operator $\LLL : BV \to BV $ acting as
\[
\LLL v= \LL_0 (v{\bf 1}_{B^c_n})=\LL(v {\bf 1}_{B^c_n}{\bf 1}_{X_0}),
\]
it is straightforward to check that
\begin{equation}\label{GK}
\mathbb{P}_n(M_n\le u_n)=\frac{1}{\al^{n-1}}\int_I \LLL^n h_0 \, dm.
\end{equation}
Roughly speaking, when $n\rightarrow \infty$, the operator $\LLL$ converges to $\LL_0$ in the spectral sense as ${\bf 1}_{B^c_n}$ becomes less and less relevant in~$\LL_0 (v{\bf 1}_{B^c_n})$. In particular, the top eigenvalue of $\LLL$ will converge to that of $\LL_0$ and this will allow us to control the asymptotic behavior of the integral on the right hand side of~\eqref{GK}. We now make these arguments rigorous by adapting the perturbative strategy put forward in~\cite{KL2, GK}. Under the following
\begin{standing} \label{def:assumptions}
Assume that $h_{-}:=\essinf_{\textrm{supp}(\Lambda)} h_0>0$ i.e. the essential
infimum is taken  with respect to~$\Lambda$. Let
\[
r_{k,n}:=
\frac{\Lambda(B_n\cap T^{-1}B^c_n\cap\cdots\cap T^{-k}B_n^c\cap
T^{-(k+1)}B_n)}{\Lambda(B_n)},
\]
where $r_{k,n}$ is the conditional probability with respect to~$\Lambda$, that we return to $B_n$ exactly after $k+1$ steps.
Assume that
\[
r_k=\lim_{n\rightarrow \infty}r_{k,n} \quad \text{exists for all } k.
%\ \sum_{k=0}^{\infty}\alpha^{-(k+1)}q_k<\infty,
\]
\end{standing}
\noindent
we will now prove that we satisfy the necessary assumptions A1--A4 of~\cite{KL2, GK}.

\noindent
{\bf A1.} The operators $\LLL$ enjoy the same
Lasota--Yorke inequalities \eqref{LY1} with the same expansion constant $\chi$
and $b$ in front of the weak norm. It is sufficient to adapt the
arguments of \cite{LM} by
replacing ${\bf 1}_{X_0}$ with ${\bf 1}_{X_0\cap B_n^c}.$

\noindent
{\bf A2}. We now compare the two operators;
here
the weak and strong Banach spaces will be again $L^1$ and $\text{BV}.$ We have:
\begin{equation}\label{FC}
\int |(\LL_0 -\LLL)v| \, dm=\int |\LL_0(v{\bf 1}_{B_n})| \, dm\le   \|v \|_{BV}
\, m(B_n\cap X_0),
\end{equation}
by expressing  $\LL_0$ in terms of $\LL$ and since the $L^{\infty}$ norm of $v$
is bounded by $\|v\|_{BV}$  in one dimensional systems,
see~\cite[Section 2.3]{GB}. Then, for the triple norm,
$ \vvvert\LL-\LLL \vvvert_1\le  m(B_n\cap X_0)$ and therefore  for $n$ large enough (see section~\ref{ssec:closeness}), we
get the following spectral properties, analogous of (\ref{QC2}), namely:
\begin{align}\label{QC3}
\LLL h_n=\lambda_n h_n, \
\LLL^* \mu_n=\lambda_n \mu_n\\
\lambda_n^{-1}\LLL g= h_n \int g \, d\mu_n + \widetilde{\mathcal{Q}}_n g,
\end{align}
where $h_n\in \text{BV}$, $\mu_n$ is a Borel measure and $\widetilde{\mathcal{Q}}_n$ a
linear operator with spectral radius less than one; moreover
%$\sup_n  \text{sp}(\widetilde{\mathcal{Q}}_n)<\text{sp}(\mathcal{Q})<1.$ \\
$\sup_n  \text{sp}(\widetilde{\mathcal{Q}}_n) <\text{sp}(\mathcal{Q})<1.$

\noindent
{\bf A3.}  Next, we need to show that
\begin{equation}\label{TTRR}
\sup \Big\{\int (\LL_0 -\LLL)v \, d\mu_0: v\in \text{BV},  \|v \|_{\text{BV}}\le1
\Big\} \times  \|\LL_0 (h_0{\bf 1}_{B_n}) \|_{\text{BV}}\le C_{\sharp} \Delta_n,
\end{equation}
where $$\Delta_n:= \int \LL_0 ({\bf 1}_{B_n}h_0) \, d\mu_0 =\alpha \Lambda(B_n)$$ and
$C_{\sharp}$ is a constant. Notice that the first term on the left hand side of (\ref{TTRR}) is the triple norm $ \vvvert\LL_0-\LLL \vvvert_{\mu_0}.$\footnote{The reader could wonder why we used two different triple norms, the first in (\ref{FC}) with respect to $m$ and the second in~\eqref{TTRR} with respect to~$\mu_0.$ The first was used to get the quasi-compactness representation for the operator $\LLL$ given in~\eqref{QC3} and we should use there the same couple of adapted functional spaces $L^1$ and $\text{BV}$ as prescribed by the main theorem in~\cite{KL2}. The second allowed us to compare the maximal eigenvalues of $\LL_0$ and $\LLL$ and it requires the eigenfunction of the dual of $\LL_0,$ which is $\mu_0$ as prescribed in~\cite{KL}. }
This is bounded by  $\alpha
\mu_0(B_n)$, as can be obtained by an argument analogous to \eqref{FC}, combined with~\eqref{DU}.\footnote{
We used here  that $\sup_I v \le v(0)+|v|_{\text{TV}}$,
where $|\cdot|_{\text{TV}}$ denotes the total variation seminorm. Since $\mu_0$ is
not atomic (see next section), we can take $v(0)=0.$ A similar estimate was used in the bound given in the proof of \cite[Lemma~7.2]{LM}.}  The second factor is bounded
by the Lasota--Yorke inequality with a constant $C_{h_0}$ depending on $h_0.$ Then by the first standing assumption  $\smash{\alpha C_{h_0}\mu_0(B_n)\le \frac{\alpha
C_{h_0}}{h_{-}}\Lambda(B_n)}$.

\noindent
{\bf A4.}
We now define the following quantity for $k\ge 0:$
\begin{equation} \begin{split} \label{QQQQ}
q_{k,n}&:=\frac{\int (\LL_0 -\LLL)\LLL^k(\LL_0 -\LLL)(h_0) \, d\mu_0}{\Delta_n}.
\end{split} \end{equation}
By the duality properties enjoyed  by the transfer operators with respect to
our standing assumption, it is easy to show that
\begin{equation}\label{QQQ}
q_{k,n}= \alpha^{k+1}r_{k,n}.
\end{equation}
We observe that by the Poincar\'e Recurrence Theorem with respect to  the invariant measure $\Lambda$, as $r_{k,n}$ is the probability that the system returns to $B_n$ in exactly $k+1$  steps, we have
\[
\sum_{k=0}^{\infty}\alpha^{-(k+1)}q_{k,n}=\sum_{k=0}^{\infty}r_{k,n} =1.
\]
We denote by $\theta$ the {\em extremal index} (EI), which will be therefore between $0$ and $1:$
\[
\theta:=1-\sum_{k=0}^{\infty}r_k.
\]
With our standing assumption, since we satisfy A1--A4, the
perturbation theorem by Keller and Liverani~\cite{KL} gives (we recall the top eigenvalue of $\LL_0,$ $\lambda_0,$
is equal to~$\alpha$)
\begin{equation}\label{AS}
\lambda_n=\alpha-\theta\ \Delta_n+o(\Delta_n)=\alpha \exp \left(-\tfrac{\theta}{\alpha}\Delta_n+o(\Delta_n)\right),\ \text{as}\  n\to\infty,
\end{equation}
or equivalently,
\begin{equation}\label{POL}
\lambda_n^n=\alpha^n\exp \left(-\tfrac{\theta}{\alpha}n\Delta_n+o(n\Delta_n) \right).
\end{equation}
We now substitute (\ref{POL}) in the right hand side of (\ref{GK}) and use (\ref{QC3})  to get
\begin{align*}
\mathbb{P}_n(M_n\le u_n) &=\frac{1}{\alpha^{n-1}}\int \lambda_n^n h_n \, dm \int
h_0 \, d\mu_n+\lambda_n^n\int \widetilde{\mathcal{Q}}_n^nh_0 \, dm \\
&= \al \exp(-\tfrac{\theta}{\alpha}n\Delta_n+o(n\Delta_n))\int h_n \, dm \int h_0
\, d\mu_n+\lambda_n^n\int \widetilde{\mathcal{Q}}_n^nh_0 \, dm. \numberthis \label{eqn}
\end{align*}
It has been proved in \cite[ Lemma 6.1]{KL}   that $\int h_0 \,
d\mu_n\rightarrow 1$ for $n\rightarrow \infty.$ Following again \cite{KL}, see
also \cite[Section 2.1]{GK}, it has been shown how to normalize
$h_n$ and $\mu_n$ in such a way that $\int h_n \, d\mu_0=1.$ But in our case we
have instead the term~$\int h_n \, dm.$ Now we observe that by \eqref{FC} and
by the perturbative theorem in \cite{KL2}, we have that $|h_n-h_0|_1\rightarrow
0$ as $n\rightarrow \infty$. Moreover
\[
\int h_0 \, dm= \frac{1}{\alpha}\int \LL _0h_0 \, dm=\frac{1}{\alpha}\int \LL
(h_0{\bf 1}_{X_0}) \, dm=\frac{1}{\alpha}\int  h_0{\bf 1}_{X_0} \,
dm=\frac{1}{\alpha}\nu(X_0)=\frac{1}{\alpha},
\]
and this term will compensate the $\al$ in the numerator in the equality above.
Note that the choice given by \eqref{BBLL}  is equivalent to
$n\Delta_n\rightarrow \alpha \tau.$ In this case $\lambda_n^n$ will be simply
bounded in $n$ and $\int |\widetilde{\mathcal{Q}}_n^n(h_0)| \, dm\le
\text{sp}(\mathcal{Q})^n  \|h_0 \|_{BV}\rightarrow 0.$ In conclusion we have
\begin{equation}\label{FF}
\lim_{n\rightarrow \infty} \mathbb{P}_n(M_n\le u_n)=e^{-\tau \theta},
\end{equation}
which is the Gumbel's law.

\section{The extremal index}
\label{sec:ei}

\subsection{Smallness of the hole}

We briefly return to the Section~\ref{ssec:closeness} to quantify
the distance between the maximal eigenvalue of $\LL$, which is $1$, and that of $\LL_0$, which is~$\alpha\le 1.$ In the previous section we described the asymptotic deviation of $\lambda_n$ from~$\alpha$ as $n\to\infty$. For the next  considerations we will compare $\alpha$ to~$1$. This is given in the \cite[formula (2.3)]{KL}, and with our notation reads as (see \eqref{CO}):
\begin{equation}\label{CA1}
1-\alpha\le \hat{C}  \vvvert\LL_0-\LL \vvvert_1\le \hat{C} m(H),
\end{equation}
where the constant $\widehat{C}$ is computed explicitly in
\cite[Section 2.1]{GK} and depends on the density~$h_0.$ We now strengthen the assumption on the ``smallness'' of the hole by requiring that $m(H)$ is such that for a fixed
$1<D<\beta=\inf_I|T'|$ it holds
\begin{equation}\label{CA2}
\alpha>\frac{D}{\beta};
\end{equation}
for instance take $m(H)\le \frac{1}{\hat{C}}(1-\frac{D}{\beta})$.
This has two interesting consequences; one will be established at the end of this section when we will compute the extremal index for periodic points. The other one states that the measure $\mu_0$, and therefore $\Lambda$, is not atomic. The proof is a straightforward adaption of \cite[Lemma~2]{HK}, where the conformal structure of $\mu_0$ is used and their ``$d$'' is replaced by our ``$D$''. Another proof of the non-atomicity of $\Lambda$ for more general holes is given in \cite[Lemma 4.3]{LM}.

\subsection{Position of the target point}
\label{ssec:PTP}

We now return to the computation of the extremal index $\theta$, which relies on the~$r_{k,n}.$ By using the fact that we restricted our considerations to the potential $\frac{1}{|T'|}$ we can easily reproduce the  arguments on the invariant set~$X_{\infty}$. These give two types of behavior according to the nature of the target point~$z$, see \cite{S1, FR1, FR2} for similar computations for different kind of dynamical systems. Recall that we write $B_n$ instead of~$B(z,e^{-u_n}).$  By recalling the definition of Lasota--Yorke maps, let $z$ be a non-periodic point and not belonging to the countable union $S$ of the preimages of the boundary points of the domains of local injectivity of $T.$ On $I\setminus S$, the maps $T^n, n\ge 1,$ are all continuous and moreover $\Lambda(I\setminus S )=1.$ Now, we fix $k$ and go to the limit for large $n$ in~(\ref{QQQ}).  By exploiting the continuity of $T^k$ and by taking $n$ large enough, all the points in $B_n$ will be around $T^k(z)$ and at a positive distance from $B_n$, so that $r_{k,n}$ is zero and no limit in $n$ is required any more.

Suppose now $z$ is a periodic point of minimal period $p;$  all the $r_{k,n}$ with $k\neq p-1$ are zero for the same reason exposed above. When $k=p-1$ any point in $B_n$ will be at a positive distance from $B_n$ when iterated $p-2$ times; this again is a consequence of continuity for large $n.$ But for $k=p-1$, $T^{k+1}(z)=T^p(z)=z;$ by taking again $n$ large enough there will be only one preimage of $T^{-p}B_n$, denoted $T_z^{-p}B_n$ intersecting $B_n.$ Since the map $T^p$ is uniformly expanding, such a preimage will be properly included in $B_n.$ We are thus led to compute
\begin{equation}\label{FD}
\frac{\Lambda(T_z^{-p}B_n)}{\Lambda(B_n)}=\frac{\int_{T_z^{-p}B_n}h_0d\mu_0}{\int_{B_n}h_0d\mu_0}.
\end{equation}
We now make an additional {\em assumption}, namely that $h_0$ is continuous at $z$; we recall that the set of discontinuity points is countable, since~$h_0\in \text{BV}.$ Since $z$ is periodic with period $p$ we have to compare the density at the numerator and at the denominator in (\ref{FD}) in two close points and both close to $z$. Therefore
\[
\frac{\Lambda(T_z^{-p}B_n)}{\Lambda(B_n)}\sim \frac{\int_{T_z^{-p}B_n}d\mu_0}{\int_{B_n}d\mu_0},
\]
and the equality will be restored in the limit of large $n$ when the previous two close points will converge to~$z$. So we are left with estimating the ratio
$\smash{\frac{\mu_0(T_z^{-p}B_n)}{\mu_0(B_n)}}$; we point out again that
$B_n=T^p(T_z^{-p}B_n)$ and that $T^p$ is one-to-one on $T_z^{-p}B_n.$
Therefore, by considering $T^p$ and iterating \eqref{CM}, we obtain
\[
\frac{\mu_0(T_z^{-p}B_n)}{\mu_0(B_n)}=\frac{\mu_0(T_z^{-p}B_n)}{\int_{T_z^{-p}B_n}\alpha^p|(T^p)'|(y)d\mu_0(y)}.
\]
Passing to the limit and exploiting again the continuity of $T^p$ at $z$, we
finally have
\[
r_{p-1}= \frac{1}{\alpha^p|(T^p)'|(z)}, \textrm{ and }
\theta=1-\frac{1}{\alpha^{p}|(T^p)'|(z)}
\]
where
$
\alpha  |T'(z)|>D>1.
$
By collecting the previous result we have proved the following:
\begin{proposition} \label{prop:final}
Let $T$ be a uniformly expanding map of the interval $I$ preserving a mixing
measure. Let us fix a small absorbing region, a hole $H\subset I$; then there
will be an absolutely continuous conditionally invariant measure $\nu,$
supported on $X_0=I\setminus H$ with density~$h_0.$ Write
$\alpha=\nu(T^{-1}X_0)$.  If the hole is small enough there will be a
probability measure $\mu_0$ supported on the surviving set $X_{\infty}$ such that
the measure $\Lambda=h_0\mu_0$ is $T$-invariant; we will assume that $h_0$ is bounded
away from zero. Having fixed the positive number $\tau$, we take the sequence $u_n$ satisfying $n\Lambda(B(z,\exp(-u_n) ))=\tau,$ where~$z\in X_{\infty}$.  Then, we take the sequence of conditional probability measures $\mathbb{P}_n(A)=\frac{\nu(A\cap X_{n-1})}{\nu(X_{n-1})},$ for $A \subset I$ measurable, and define the random variable $M_n (x):=\max\{\phi(x), \cdots,
\phi(T^{n-1}x)\},$ where $\phi(x)=-\log|x-z|.$ Moreover we will suppose that all the iterates $T^n, n\ge 1$ are continuous at $z$ and also that $h_0$ is continuous
at $z$ when the latter is a  periodic point. Then we have:
\begin{itemize}
\item If $z$ is not a periodic point:
\[
\mathbb{P}_n(M_n\le u_n)\rightarrow e^{-\tau}.
\]
\item If $z$ is a periodic point of minimal period $p$, then
\[
\mathbb{P}_n(M_n\le u_n)\rightarrow e^{-\tau\theta},
\]
where the extremal index $\theta$ is given by:
\[
\theta=1-\frac{1}{\alpha^{p}|(T^p)'|(z)}
\]
\end{itemize}
\end{proposition}

Note that in literature, the {\em escape rate} $\eta$  for our
open system is usually defined as $ \eta=-\log \alpha$ thus we can
see the extremal index as
\[
\theta= 1-\frac{1}{e^{-p
\eta}|(T^p)'|(z)}.
\]

\begin{remark} \label{rem:extension1}
We presented here the simplest possible case. However,
starting again from the transfer operator \eqref{transferop}, it could be possible to perform the same analysis with a generic potential,  adapting the construction of the spaces to handle different weights. As a starting point, \cite{LM} contains elements to treat conditional measures in such situation.
\end{remark}

\begin{remark} \label{rem:extension2}
In light of \cite{BDT,DT}, it would be interesting to construct a statement analogous to our main Proposition~\ref{prop:final}, when either the hole is not of a given size or either the dynamics generated is mixing at a subexponential rate.
\end{remark}

\begin{remark} \label{rem:billiard}
An analogous billiard statement, following \cite{D}, could be constructed from
the above provided there is enough hyperbolicity to beat the complexity growth.
In a nutshell, given a billiard, one can consider the Poincar\'e map given by
the collision with the scatterers. One has then the freedom to choose absorbing
scatterers and target scatterers as long as the absorbing part is not too wide.
\end{remark}

\begin{remark}  \label{rem:linearresponse}
 As the approach to study the extremal index is perturbative in nature, should
not come as a surprise that one could consider a one-parameter family of maps
$T_\varepsilon$ which are small perturbations of~$T$. It could be possible,
following some of the techniques of \cite{S1,GG}, to establish the
behaviour of the extremal index with respect to deterministic perturbations or noisy perturbations in our framework of targets and holes.
\end{remark}

\subsection{On the choice of the boundary sequence}
\label{ssec:BS}

Let us now comment on \eqref{BBLL}, i.e., the scaling behavior~$n\Lambda(B(z, e^{-u_n}))\rightarrow \tau$. Since as we already argued,  the
measure $\Lambda$ is not atomic, it varies continuously with the radius of the
ball; therefore for any fixed $\tau$ and $n$ we could choose $u_n$ so that
\begin{equation}\label{BLS}
\Lambda(B(z,e^{-u_n}))=\frac{\tau}{n}.
\end{equation}
Unluckily,  the measure
$\Lambda$ is often not computationally accessible; however we can use the following approximation scheme to construct a sequence of $u_n$ which still satisfies \eqref{BBLL}. Let
\[
d_n(z):= \frac{\log \Lambda(B(z,e^{-u_n}))}{\log e^{-u_n}}.
\]
Since the density $h_0$ is bounded away from zero by the standing assumptions, for $\delta$ arbitrarily small and  $n$ large enough we have that
\[
d_n(z)\ge \frac{\log \mu_0(B(z,e^{-u_n}))}{\log e^{-u_n}}-\delta.
\]
By \cite[Theorem B]{LM}, whenever the map $T$ has large
images and large images with respect to the hole $H$ (see the
discussion before \cite[Theorem B]{LM}), then for all
$z\in X_{\infty}$, there exists $t_0>0$ such that 
\[
\liminf_{n\rightarrow \infty}\frac{\log \mu_0(B(z,e^{-u_n}))}{\log e^{-u_n}}\ge t_0
\]
and the Hausdorff dimension of the surviving set $HD(X_{\infty})$ verifies
\[
HD(X_{\infty})\ge t_0.
\]
Therefore if we fix again $\delta$ and take correspondingly $n$ large enough we have  that
$
d_n(z)\ge t_0-\delta-\delta\ge t_0-2\delta
$
which implies
$
\Lambda(B(z,e^{-u_n}))\le e^{-u_n (t_0-2\delta)},
$
and, together with \eqref{BLS}, finally
$
 \tau \le n e^{-u_n (t_0-2\delta)}.
$
In other words, $\smash{u_n\le -\frac{\log \tau}{t_0-2\delta}+\frac{\log n}{t_0-2\delta}}$, which can also be written as
\begin{equation}\label{FDD}
\sup_n\big\{u_n-\tfrac{\log n}{t_0} \big\}\le -\frac{\log \tau}{t_0}.
\end{equation}
In the computational approach to extreme value theory, the boundary level
$u_n$ are chosen with the help of an affine
function (see~\cite{Led}):
\[
u_n=\frac{\log\tau^{-1}}{a_n}+b_n.
\]
The sequences $a_n$   and $b_n$ can be obtained with the help of the Generalized Extreme Value (GEV) distribution in order to fit Gumbel's law. The inequality \eqref{FDD} suggests that for $n$ large  $a_n\sim t_0$ and $b_n\sim \frac{\log n}{t_0}$, therefore we could attain a lower bound for the Hausdorff dimension of the surviving set. We defer, for instance, to \cite{v} to show  how
to use the GEV distribution to estimate the sequences $a_n, b_n$, and we will show in future studies how to use such estimates to approach~$HD(X_{\infty}).$

\section{How far are we from the survining set? The degenerate limit}
\label{hf}

We noted several times that the support of $\mu_0$ is the surviving set~$X_{\infty}$. This means that if we pick the open ball $B_n=B(z, e^{-u_n})$ centered in a point $z\notin X_{\infty}$  or even in the hole, then when the radius of the ball is sufficiently small, we have~$\mu_0(B_n)=0$, since~$X_{\infty}$ is a closed set. This immediately implies by the argument similar to that we used in (\ref{CA1}) that
$$
|\lambda_n-\alpha|\le \text{const} \times  \vvvert\LL_0-\LLL \vvvert_{\mu_0}\le  \text{const}\times \alpha \mu_0(B_n)=0.
$$
 The fact that the perturbed eigenvalue could become equal to the unperturbed one for a finite size of the perturbation, is already a part of \cite[Theorem 2.1 ]{KL} and is also detailed in~\cite[Footnote (3)]{GK}. Therefore, if we call $\hat{n}$ the first $n$ for which $B_n\cap X_{\infty}=\emptyset,$ for any $n\ge \hat{n}$ we have that
\[
\mathbb{P}_{n}(M_{n}\le u_n)=\alpha \int h_{n} \, dm \int h_0
\, d\mu_{n}+\alpha^{n}\int \widetilde{\mathcal{Q}}_{n}^{n}h_0 \, dm.
\]
As explained above, for $n\rightarrow \infty$ it holds $\int h_0
\, d\mu_n\rightarrow 1$, $\int h_n \ dm\rightarrow \alpha^{-1}$, and $\int \widetilde{\mathcal{Q}}_{n}^{n}h_0 \, dm\rightarrow 0$, we thus have that
\begin{equation}\label{NIB}
\mathbb{P}_{n}(M_{n}\le u_n)\rightarrow 1, \ n\rightarrow \infty.
\end{equation}
Trivially, \eqref{NIB} states that if the target point is off the surviving set, then the trajectories will not be able to approach it arbitrary close.
This result has two interesting consequences for applications, in particular the second one will provide a full description of the extreme value distribution (EVD) for any choice of the target set.

First, we observe that the limit (\ref{NIB}) holds for any sequence $u_n$ going to infinity, and for simplicity we now put~$u_n=\log n$. Then we could reasonably argue that for the smallest $\hat{n}$ for which
\[
\mathbb{P}_{\hat{n}}(M_{\hat{n}}\le \log\hat{n})\sim 1,
\]
then
\[
\text{dist}(z, X_{\infty})\sim \frac{1}{\hat{n}}.
\]

Second, let us return to the statement of our main Proposition~\ref{prop:final}. Whenever
we take the point $z\in X_{\infty}$  and by a suitable choice of the sequence $u_n$ as we explained in Section~\ref{ssec:BS}, we get a non-degenerate limit for our EVD, in particular different from~$1.$ Instead, if we pick the point $z$ outside the surviving set and no matter what the sequence $u_n$ is,
provided it goes to infinity, we get a degenerate limit equal to one for the EVD.

\section{Aknowledgements}
P.G. acknowledges the support of the Centro di Ricerca
Matematica Ennio de Giorgi and of UniCredit Bank R\&D group for financial
support through the ``Dynamics and Information Theory Institute'' at the
Scuola Normale Superiore. P.K. is supported by the Deutsche
Forschungsgemeinschaft (DFG) through the CRC 1114 ``Scaling Cascades in Complex
Systems'', project A01; he  thanks also the  {\em Centro de Giorgi} in Pisa
where this work was initiated.
S.V. thanks the Laboratoire International Associ\'e LIA LYSM, the INdAM
(Italy), the UMI-CNRS 3483, Laboratoire Fibonacci (Pisa) where this work has
been completed under a CNRS delegation and the {\em Centro de Giorgi} in Pisa
for various supports.


\begin{thebibliography}{10}
\bibitem{S1} H. Aytac, J. Freitas, S. Vaienti, Laws of rare events for
deterministic and random dynamical systems, {\em Trans. Amer Math. Soc.}, 367
(2015), 8229-8278.
\bibitem{BV} W. Bahsoun, S. Vaienti, Escape rates formulae and metastability for randomly perturbed maps, {\em Nonlinearity}, 26 (2013) 1415--1438
\bibitem{GB} A. Boyarsky,  P. F. G{\'o}ra, {\em Laws of Chaos: Invariant Measures
and Dynamical Systems in One Dimension},  (Birkh\"auser)  1997
 doi:10.1007/978-1-4612-2024-4
\bibitem{BDT} H. Bruin, M. F. Demers , M. Todd, Hitting and escaping
statistics: mixing, targets and holes {\em Advances in Mathematics} 328 (2016)
doi:10.1016/j.aim.2017.12.020
\bibitem{S3} Th. Caby, D. Faranda, G. Mantica, S. Vaienti, P. Yiou,
Generalized dimensions, large deviations and the distribution of rare events,  {\em PHYSICA D}, 400, 132-143, (2019).
\bibitem{S4}Th. Caby, D. Faranda, S. Vaienti, P. Yiou,  On the
computation of the extremal index for time series, submitted,
https://arxiv.org/pdf/1904.04936.pdf
\bibitem{D} M.F. Demers, Dispersing Billiards with Small Holes,
 in "Ergodic Theory, Open Dynamics, and Coherent Structures" Springer New
York 2014 ,  doi:10.1007/978-1-4939-0419-8\_8 , isbn=978-1-4939-0419-8

\bibitem{DT} M.F. Demers, M. Todd , Slow and Fast Escape for Open Intermittent
Maps, {\em Communications in Mathematical Physics } 351 (2017)
doi:10.1007/s00220-017-2829-6
\bibitem{S2} D. Faranda, H, Ghoudi, P. Guiraud, S. Vaienti, Extreme value
theory for synchronization of coupled map lattices, {\em Nonlinearity},  31, 7,
3326-3358 (2018).
\bibitem{v} D. Faranda, V. Lucarini, G. Turchetti, S. Vaienti, Extreme value distributions for singular measures, {\em Chaos}, 22, 023135, (2012).
\bibitem{FR1}A.C.M. Freitas, J.M. Freitas, M. Todd, Hitting Time Statistics and Extreme Value Theory, {\em Probab. Theory Related Fields}, 147, no. 3, (2010), 675-710.
    \bibitem{FR2} A.C.M. Freitas, J.M. Freitas, M. Todd, Extremal Index, Hitting Time Statistics and periodicity, {\em Adv. Math}., 231, no. 5, (2012), 2626-2665.
\bibitem{FFT} A. Freitas, J. Freitas and M.J. Todd,  Speed of convergence for
laws of rare events and escape rates, {\em Stochastic Processes and their
Applications" } 125 (2015) doi:"10.1016/j.spa.2014.11.011".

\bibitem{GG} G. Galatolo, P. Giulietti,  A linear response for dynamical systems
with additive noise, {\em Nonlinearity}  32 (2019)   doi:
10.1088/1361-6544/ab0c2e
\bibitem{HK} F. Hofbauer, G. Keller, Ergodic properties of invariant measures for piecewise monotonic transformations, {\em Math. Zeit}, 180, 119-140, (1982).

\bibitem{GK}G. Keller, Rare events, exponential hitting times and extremal indices via spectral perturbation, {\em Dyn. Syst.},
27 (2012) 11-27.
\bibitem{KL} G. Keller and C. Liverani, Rare events, escape rates and quasistationarity: some exact formulae, {\em J. Stat.
Phys.},  135 (2009), 519-534.
\bibitem{KL2} G. Keller, C. Liverani, Stability of the spectrum for transfer operators.
    {\em Annali della Scuola Normale Superiore di Pisa, Classe di Scienze}, (4) Vol. XXVIII, pp. 141--152 (1999).
    \bibitem{Led} G.M.R. Leadbette, H. Rootzen, {\em Extremes and Related Properties of Random Sequences and
Processes}, Springer Series in Statistics, New York: Springer, (1983).
\bibitem{LM} C. Liverani, V. Maume-Deschamps, Lasota--Yorke maps with holes: conditionally invariant probability measures and invariant probability measures on the survivor set, {\em Annales de l'Institut Henri Poincar\'e Probability and Statistics}, 39 (3), 385-412 (2003).
\bibitem{S5} Lucarini  V, Faranda  D, Freitas  A  M, Freitas  J  M,
Holland  M, Kuna  T, Nicol  M, Todd  M and
Vaienti S, {\em Extremes and Recurrence in Dynamical Systems} (New York: Wiley),
2016





    




\end{thebibliography}
\end{document}